
\magnification=\magstep1
\input amstex
\documentstyle{amsppt}



\define\Spec{\operatorname{Spec}}
\define\id{\operatorname{id}}

\define\res{\operatorname{res}}

\define\A {\bold A}
\topmatter
\title\nofrills Symmetric tensors with  applications to Hilbert schemes\\
\endtitle 
\author Roy Mikael Skjelnes\endauthor
\affil Department of Mathematics,  KTH \endaffil
\address KTH, S-100 44 STOCKHOLM, Sweden\endaddress
\email skjelnes\@math.kth.se \endemail 
\subjclass{14C05, 14D22}\endsubjclass

\abstract 
Let $A[X]_U$ be a fraction ring of the polynomial ring $A[X]$ in the
variable $X$ over a commutative ring $A$. We show that the Hilbert
functor ${\Cal Hilb}^n_{A[X]_U/A}$ is represented by an affine scheme
$\text{Symm}^n_A(A[X]_U)$ given as the ring of symmetric tensors of
$\otimes_A^nA[X]_U$. The universal family is given as
$\text{Symm}^{n-1}_A(A[X]_U)\times_A \Spec(A[X]_U)$.
\endabstract
\endtopmatter
\document

\subhead \S 1. - Introduction\endsubhead 

The Hilbert functor ${\Cal Hilb}^n_{{\Cal O}_{V,x}}$ parameterizing
closed subschemes of a
variety $V$, having finite length $n$ and support in a fixed  point $x
\in V$, has been studied by several authors in the last
decades (\cite{2}, \cite{3}, \cite{4}, \cite{6}, \cite{7}, \cite{8}  and \cite{13}). The main
interest have been on the classification of the set of $k$-rational
points of ${\Cal Hilb}^n_{{\Cal O}_{V,x}}$.  The scheme structure on the parameter schemes was apparently
neglected until \cite{13}, where the Hilbert functor ${\Cal Hilb}^n_{k[X]_{(X)}}$ parameterizing
subschemes of the line of finite length $n$ and support in the origin, was
described. 

It was shown in \cite{13} that the scheme representing the functor ${\Cal
  Hilb}^n_{k[X]_{(X)}}$ is of dimension $n$. Taking into account that the set of $k$-rational points of ${\Cal
  Hilb}^n_{k[X]_{(X)}}$ is trivial, the scheme structure on the
  parameter scheme was a surprise. The Hilbert scheme parameterizing
  subschemes of the line, having finite length and support in one
  point is not an algebraic scheme, that is a scheme which is not of finite type over
  the base.

The motivation behind the present paper arose from the desire to
better understand the techniques introduced in \cite{13}. Instead of
only consider finite length subschemes of the line with support in a
fixed point, we allow the subschemes to have support in any given set. 

Let $A[X]_U$ be the fraction ring of the polynomial ring $A[X]$ in the
variable $X$ over a base ring $A$, with respect to a multiplicatively
closed subset $U\subseteq A[X]$. We study the contravariant functor
${\Cal Hilb}^n_{A[X]_U/A}$ from the category of $A$-schemes to sets,
which sends an $A$-scheme $T$ to the set of closed subschemes
$Z\subseteq T\times_A \Spec(A[X]_U)$ such that the projection map $p :
Z \to T$ is flat and where the global sections of the fiber
$p^{-1}(t)$ is of rank $n$ for all points $t\in T$.

Our main result is Theorem (8.2), where we show that ${\Cal Hilb}^n_{A[X]_U/A}$ is
represented by the $n$-fold symmetric product
$\text{Symm}_A^n(A[X]_U)$ of $\Spec(A[X]_U)$. Thus we obtain a result which is similar to
what is known for Hilbert schemes of points on smooth projective curves
\cite{5},
and similar with the results of B. Iversen about the $n$-fold sections for smooth families of
curves \cite{9}.

When studying the functor ${\Cal Hilb}^n_{A[X]_U/A}$ the key problem is to determine those monic polynomials
$F(X)$ in $A[X]$ such that the fraction map $A[X]/(F(X)) \to
A[X]_U/(F(X))$ is an isomorphism. A problem which is solved by the use
of the symmetric operators \cite{13}  of the polynomial ring $A[X]$ associated to
$F(X)$. We define, Section (2.2), an $A$-algebra homomorphism $u_F$ from the ring of
symmetric  tensors of $\otimes_A^nA[X]$
to $A$, which has the following property. For every $f(X)$ in $A[X]$ 
the residue class of $f(X)$ modulo the ideal $(F(X))$ gives
by multiplication an $A$-linear endomorphism $\mu_F(f)$ on the free
$A$-module $A[X]/(F(X))$. In Theorem (2.4) it is shown that
$u_F(f(X)\otimes \dots \otimes f(X))=\text{det}(\mu_F(f))$. 

Theorem (2.4) is the technical heart of the present paper and gives a
nice relation between the symmetric operators  with
the coefficients of the characteristic polynomial of $\mu_F(f)$. A relation of the
homomorphism $u_F$ to resultants is touched upon in Section (3).

By studying the properties of such homomorphisms $u_F$ we classify in Section (4), the set of ideals $I\subseteq
A[X]_U$ such that the residue ring $A[X]_U/I$ is free as an
$A$-module. 

We denote with $\otimes_A^{(n)}A[X]$ the ring of invariant tensors of
$\otimes_A^nA[X]$ under the standard action of the symmetric group in
$n$-letters. In Section (5) we show that the functor parameterizing
ideals of $A[X]_U$ such that the residue rings are free of rank $n$ as $A$-modules is
represented by a fraction ring of $\otimes_A^{(n)}A[X]$. A fraction
ring which we show in Section (6) is isomorphic to the ring of
symmetric tensors $\otimes_A^{(n)}A[X]_U$.

Then finally we summarize the accumulated results in the main Theorem
which states that the $n$-fold symmetric product
$\text{Symm}^n_A(A[X]_U)$ represents the functor ${\Cal
  Hilb}^n_{A[X]_U/A}$, and where the universal family is given as
$\text{Symm}^{n-1}_A(A[X]_U)\times_A \Spec(A[X]_U)$.

As an  application of our result we have  that the Hilbert scheme of $n$-points on the
affine line $\A^1_A$ over $\Spec{A}$ is represented by the affine
$n$-space $\A^n_A$. Here $A$ is any commutative unitary
ring. Furthermore we have that ${\Cal Hilb}^1_{A[X]_U/A}$, the Hilbert
functor parameterizing 1 point on $\Spec(A[X]_U)$ is represented by
the scheme $\text{Spec}(A[X]_U)$, \cite{10}. In the case when the
base ring $A=k$ is a field, we have that the ring of 
symmetric tensors of $\otimes_k^nk[X]_{U}$ parameterizes the closed
subschemes of the line $\text{Spec}(k[X])$ having length $n$ and
support in the subset of the line corresponding to the prime ideals   $\{P \subset k[X] \mid P\cap U
=\emptyset \}$. In particular
we recover the situation considered in  \cite{13} where $k[X]_U=k[X]_{(X)}$ is the local ring
of the maximal ideal $(X)\subseteq k[X]$. A more precise discussion of
some applications of Theorem (8.2) is given in the end of Section (8).

I thank D. Laksov for the help and assistance I have received during
our many discussion on the present subject. I thank S. A. Str\o mme for his
comments and remarks.

\subhead \S 2. - Symmetric operators on the polynomial ring \endsubhead

Given a monic polynomial $F(X)$ in $A[X]$ which has positive degree
$n$. Then the residue ring $A[X]/(F(X))$ is a free $A$-module of rank
$n$. Any element $f(X)$ in $A[X]$ gives by multiplication by the
residue  class
of $f(X)$ modulo the ideal  $(F(X))$, an $A$-linear endomorphism $\mu_{F}(f)$ on
$A[X]/(F(X))$. The goal of Section (2) is the
construction of  an $A$-algebra homomorphism 
$ u_F : \otimes_A^{(n)}A[X] \to A$, where $\otimes_A^{(n)}A[X]$ is the
ring of symmetric functions,
with the property that for any $f(X)$ in $A[X]$ we have that
$$ u_F(f(X)\otimes \dots \otimes f(X))=\det (\mu_F(f)).$$

\definition{2.1. The symmetric operators} We recall the symmetric
operators which were introduced in \cite{13}. Let $A[X]$ denote the ring of polynomials in the variable $X$
over a commutative ring $A$.  We write $\otimes_A^nA[X]$ for the tensor algebra
$A[X]\otimes_A \cdots \otimes _A A[X]$ ($n$-copies of $A[X]$). To
every element $f(X)$ in $A[X]$ we let $f(X_i)=1\otimes \cdots \otimes
1 \otimes f(X) \otimes 1\otimes \cdots \otimes 1$ in
$\otimes_A^nA[X]$, where the $f(X)$ occurs at the $i$'th place. 
We identify $X=1\otimes\dots 1\otimes X$
in $\otimes_A^{n+1}A[X]$, and we consider $\otimes_A^{n+1}A[X]$ as the
polynomial ring in the variable $X$ over $\otimes_A^nA[X]$.

To each element $f(X)$ in $A[X]$ and for every positive integer $n$,
we define the symmetric tensors $s_{1,n}(f(X)), \ldots ,
s_{n,n}(f(X))$ by the following identity in
$\otimes_A^{n+1}A[X]$.
$$\aligned 
\Delta_{n,f(X)}(X)&= \prod_{i=1}^n (X-f(X_i)) \\
 &= X^n-s_{1,n}(f(X))X^{n-1}+\dots +(-1)^ns_{n,n}(f(X)). 
\endaligned
\tag{2.1.1}$$
We
have that $s_{1,n}(X), \ldots ,s_{n,n}(X)$ are the elementary
symmetric functions in $X_1, \ldots , X_n$. We denote the ring of
symmetric tensors of $\otimes_A^nA[X]$, which is the polynomial ring in the variables $s_{1,n}(X), \ldots
,s_{n,n}(X)$ over $A$, by
$\otimes_A^{(n)}A[X]$.

The element
$\Delta_{n,f(X)}(X)$ is a polynomial in the variable $X$,
having coefficients in the ring of symmetric functions
$\otimes_A^{(n)}A[X]$. Since $\Delta_{n,f(X)}(X)$ is a monic polynomial
of degree $n$, we have that the residue ring 
$$ V_{n,f(X)}=\big( \otimes_A^{(n)}A[X] \otimes_A
A[X]\big) /(\Delta_{n,f(X)}(X))\tag{2.1.2}$$
is a free $\otimes_A^{(n)}A[X]$-module of rank $n$.
\enddefinition

\definition{2.2.  The homomorphism $u_F$} Let
$F(X)=X^n-u_1X^{n-1}+\dots +(-1)^nu_n$ be a monic polynomial in the
polynomial ring
$A[X]$. There is a unique $A$-algebra homomorphism
$$ u_F : \otimes _A^{(n)}A[X] \to A$$
determined by $u_F(s_{i,n}(X))=u_i$. The homomorphism $u_F$ gives
$A[X]$ an $\otimes_A^{(n)}A[X]$-module structure and we have a natural
identification 
$$V_{n,X}\otimes_{\otimes_A^{(n)}A[X]}A\cong A[X]/(F(X)).$$
\enddefinition

\proclaim{Theorem 2.3} Let $M$ be a quadratic matrix, having
coefficients in a commutative ring $A$. Assume that the characteristic
polynomial $P_M(X)=\prod_{i=1}^n(X-a_i)$ of $M$ splits into linear
factors over $A$. Then we have for any polynomial $f(X)$ in $A[X]$
that the matrix $f(M)$ has characteristic polynomial
$P_{f(M)}(X)=\prod^n_{i=1}(X-f(a_i))$.
\endproclaim

\demo{Proof} The result is well known when $A$ is a field. A proof of
the Spectral Theorem over general commutative rings is found in
\cite{14}.
\enddemo

\proclaim{Theorem 2.4} Let $F(X)$ be a non-constant, monic polynomial in
$A[X]$. Denote with $n$  the degree of $F(X)$. For any element $f(X)$ in
$A[X]$ we let $\mu_F(f)$ be the $A$-linear endomorphism on $A[X]/(F(X))$
given as multiplication by the residue class of $f(X)$ modulo the
ideal $(F(X))$. We have that the
characteristic polynomial of $\mu_F(f)$ is
$$ X^n-u_F(s_{1,n}(f(X)))X^{n-1}+\dots +(-1)^nu_F(s_{n,n}(f(X))).$$
In particular we have that $u_F(f(X)\otimes \dots \otimes f(X))$ is
the determinant of $\mu_F(f)$.
\endproclaim

\demo{Proof}  For any $f(X)$ in $A[X]$
we let $\mu(f)$ be the $\otimes_A^{(n)}A[X]$-linear endomorphism on
$V_{n,X}=\otimes_A^{(n)}A[X] \otimes_AA[X]/(\Delta_{n,X}(X))$ given as
multiplication by the residue class of $f(X)$ in
$V_{n,X}$. Let $u_F : \otimes_A^{(n)}A[X] \to A$ be the $A$-algebra
homomorphism determined by $F(X)$ in $A[X]$. We have that the induced
$A$-linear endomorphism $\mu(f)\otimes \id$ on
$V_{n,X}\otimes_{\otimes_A^{(n)}A[X]}A\cong A[X]/(F(X))$ is
$\mu_F(f)$. Hence to prove our Theorem it suffices to show that the
endomorphism $\mu (f)$ has characteristic polynomial
$$ X^n-s_{1,n}(f(X))X^{n-1}+\dots +(-1)^ns_{n,n}(f(X)).\tag{2.4.1} $$
To show that $\mu (f)$ has characteristic polynomial (2.4.1) I claim
that it is
sufficient to show  that $\mu (X)$ has characteristic
polynomial $\Delta_{n,X}(X)$. Indeed, we have that $\otimes_A^{(n)}A[X]$ is a subring of
$\otimes_A^nA[X]$. Hence by ring extension, we consider $\mu (f)$ as an
$\otimes_A^{n}A[X]$-linear endomorphism on 
$$V_{n,X}\otimes_{\otimes^{(n)}_AA[X]}\otimes_{A}^nA[X]\cong
\otimes_A^{n+1}A[X]/(\Delta_{n,X}(X)). \tag{2.4.2} $$
We have that that $\Delta_{n,X}(X)=\prod_{i=1}^n(X-X_i)$ splits into
linear factors over $\otimes_A^nA[X]$. Thus if $\mu (X)$ has
characteristic polynomial $\Delta_{n,X}(X)$ then it follows from
Theorem (2.3) that $\mu (f)$ has characteristic polynomial
$\Delta_{n,f(X)}(X)=\prod_{i=1}^n(X-f(X_i))$. We have that
$\Delta_{n,f(X)}(X)$ written out in terms of symmetric functions is
(2.4.1). Thus what remains is to show that the endomorphism $\mu (X)$ on
$\otimes_A^{n+1}A[X]/(\Delta_{n,X}(X)$ has characteristic polynomial
$\Delta_{n,X}(X)$.

Let $x^i$ be the residue class of $X^i$ modulo the
ideal $(\Delta_{n,X}(X))$ in $\otimes_A^{n+1}A[X]$. We have that $1, x, \ldots ,
x^{n-1}$ form a $\otimes_A^{n}A[X]$ basis for (2.4.2). The matrix
$M$ representing the endomorphism $\mu(X)$, with respect to
the given basis is easy to describe and is called the {\it companion
  matrix} of $\Delta_{n,X}(X)$. In general, if
$F(X)=X^n-u_1X^{n-1}-\cdots -u_n$ is a monic polynomial, then the 
  companion matrix of $F(X)$ is the matrix
$$M_F =\left( \matrix
 0 & 0 & \dots &0 & u_n \\
 1 & 0 &       &\vdots       & u_{n-1} \\
 0 & 1 & \ddots & \vdots     & \vdots  \\
\vdots &  & \ddots & 0 & \vdots \\
0 &\cdots &\cdots &1 & u_1     \endmatrix \right).$$
Note that the matrix obtained by deleting first row
and first column of $M_F$, is the companion matrix of
$G(X)=X^{n-1}-u_1X^{n-2}-\cdots -u_{n-1}$. It follows readily by
induction on the size $n$ of $M_F$, that the determinant $\det
(XI-M_F)=F(X)$. Thus we get that the matrix $M$ representing the
endomorphism $\mu(X)$ with respect to the basis $1, x,
\ldots , x^{n-1}$ has characteristic polynomial $\Delta_{n,X}(X)$. We have proven the Theorem.
\enddemo

\subhead \S 3. - The Norm function and resultants\endsubhead

We will use the notation from the preceding sections. Let $\varphi :
A \to K$ be an ring homomorphism. For any
$f(X)=a_mX^m+\dots +a_m$ in $A[X]$, we write $f^{\varphi}(X)=\varphi
(a_m)X^m+\dots +\varphi (a_m)$ in
$K\otimes_A A[X]=K[X]$. Thus $f^{\varphi}(X)$ is the reduction of
$f(X)$ modulo the kernel of $\varphi$.

\definition{3.1. Definition} Given a monic polynomial $F(X)$ in
$A[X]$ which  has positive degree $n$. We define the {\it Norm function} $N_F : A[X] \to A$
with respect to $F(X)$, by sending $f(X)$ in
$A[X]$ to 
$$ N_F(f(X)):=u_F(f(X)\otimes \cdots \otimes
f(X)), \tag{3.1.1}$$
where $u_F : \otimes_A^{(n)}A[X] \to A$ is the $A$-algebra
homomorphism (2.2) determined by $F(X)$. Let $\mu_F(f)$ be the
$A$-linear endomorphism on $A[X]/(F(X))$ given as multiplication by
the residue class of $f(X)$ modulo the ideal $(F(X))$. By Theorem (2.4) we have that
$N_F(f)$ is the determinant of the endomorphism $\mu_F(f)$. We say that $N_F(f(X))$ is the {\it
  norm} of $f(X)$ with respect to $F(X)$. 
\enddefinition

Let $\varphi : A \to K$ be an $A$-algebra homomorphism. If $F(X)$ is
a monic polynomial of degree $n$ in $A[X]$, then $F^{\varphi}(X)$ is a
monic polynomial of degree $n$ in $K[X]$. Therefore $F^{\varphi}(X)$
determines a $K$-algebra homomorphism $u_{F^{\varphi}} : \otimes_K^{(n)}
K[X] \to K$. We have the following relation between the norm function
with respect to $F(X)$ and the norm function with respect to
$F^{\varphi}(X)$.

\proclaim{Lemma 3.2} Let $\varphi : A \to K$ be an $A$-algebra
homomorphism. Let $F(X)$ in $A[X]$ be a monic polynomial of positive degree $n$. Then for all $f(X)$ in $A[X]$
we have for each $i=1, \ldots , n$ that $\varphi \circ
u_F(s_{i,n}(f(X)))=u_{F^{\varphi}}(s_{i,n}(f^{\varphi}(X)))$.
In particular we have that $\varphi (N_F(f(X)))=N_{F^{\varphi}}(f^{\varphi}(X))$.
\endproclaim
\demo{Proof} The $A$-algebra homomorphism $\varphi : A \to K$ induces
by base change a homomorphism $\hat \varphi :\otimes_A^{(n)}A[X] \to
\otimes_K^{(n)}K[X]$, which maps $s_{i,n}(f(X))$ to
$s_{i,n}(f^{\varphi}(X))$ for each $i=1, \ldots ,n$ and for each
$f(X)$ in $A[X]$. Furthermore it is clear that $\varphi \circ
u_F=u_{F^{\varphi}}\circ \hat
\varphi$. where $u_{F^{\varphi}}$ is the
$K$-algebra homomorphism $u_{F^{\varphi}} : \otimes_K^{(n)}K[X] \to K$
determined by the monic polynomial $F^{\varphi}(X) \in K[X]$. We have proven the Lemma.
\enddemo

\proclaim{Proposition 3.3} Let $P(X)$ and
$Q(X)$ be two monic polynomials in $A[X]$. Let $p$ be the degree of $P(X)$ and let $q$ be the
degree of $Q(X)$. Assume that both $p$ and $q$ are positive. Then we have that
$N_P(Q(X))=(-1)^{pq}N_{Q}(P(X))$.
\endproclaim

\demo{Proof} Let $X_i=1\otimes \dots \otimes 1 \otimes X \otimes 1
\dots \otimes 1$ in the ring $\otimes_A^{p+q}A[X]$, where the $X$
occurs at the $i$'th place. We consider the
following product
$$ \res (p,q)=\prod_{i=1}^p \prod_{j=1}^q (X_i-X_{p+j}). \tag{3.3.1}$$
It follows from (3.3.1) that
the product $\res (p,q)$ is symmetric in  $X_1, \ldots , X_p$ and
symmetric in  $X_{p+1}, \ldots ,X_{p+q}$. Thus $\res(p,q)$ is an
element of $\big(\otimes_A^{(p)}A[X]\big)\otimes_A \big(\otimes_{A}^{(q)}A[X]\big)$. The monic polynomial $P(X)$ in $A[X]$ determines an $A$-algebra
homomorphism (2.2) $u_P :\otimes_A^{(p)}A[X] \rightarrow A$. We let $\hat u_P : \big(\otimes_A^{(p)}A[X]\big)
\otimes_A \big(\otimes_A^{(q)}A[X]\big) \to \otimes_A^{(q)}A[X]$ be the
induced map. Similarily we let $\hat u_Q : \big(\otimes_A^{(p)}A[X]\big)
\otimes_A \big(\otimes_A^{(q)}A[X]\big) \to \otimes_A^{(p)}A[X]$ be the
map induced by $u_Q : \otimes_A^{(q)}A[X] \to A$. Clearly we have
that $u_P\circ \hat u_Q=u_Q\circ \hat u_P$.

For fixed $i$ we have that
$$ \prod_{j=1}^q (X_i-X_{p+j})=X^q_i-s_{1,q}(X)X_i^{q-1}+\dots
+(-1)^qs_{q,q}(X). \tag{3.3.2} $$
It follows from (3.3.2) and the definition of the homomorphism $u_Q$ that $\hat u_Q$ maps $\text{res}(p,q)$ to $\prod_
{i=1}^p(Q(X_i))=Q(X)\otimes \dots \otimes Q(X)$ in  $\otimes_A^{(p)}A[X]$. Similarily we get that $\res (q,p)=(-1)^{pq}\res (p,q)$ is
mapped to $\prod_{j=1}^qP(X_{j})$ in $\otimes_A^{(q)}A[X]$ by $\hat u_P$. Thus
we have that
$$ \aligned
N_P(Q(X)) &= u_P(Q(X)\otimes \cdots \otimes Q(X))=u_P(\hat u_Q(\res
(p,q)) \\
 &=u_Q(\hat u_P(\res (p,q))) =u_Q(\hat u_P((-1)^{pq}\res(q,p))) \\
 &=(-1)^{pq}u_Q(P(X) \otimes \cdots \otimes P(X))=
 (-1)^{pq}N_Q(P(X)), \endaligned $$
proving our claim.   
\enddemo 
 
\definition{Remark} When $A$ is a field we have that  Proposition (3.3) is
a well-known formula for resultants.
\enddefinition 
\definition{Remark} The trick which we use in the proof of the Proposition
is to consider the product $\operatorname{res}(p,q)$ (3.3.1), which we found in \cite{11} (Chapter IV, \S 8, the proof of
Proposition (8.3), pp. 202-203).  
\enddefinition

\subhead \S 4. - Residues of fraction rings\endsubhead

In this section we investigate and describe ideals $I$ in fraction
rings $A[X]_U$, such that the residue ring $A[X]_U/I$ is a free
$A$-module of rank $n$. 

The key point is
Theorem (4.2) below which generalizes the Main theorem of
\cite{13} (Theorem (2.3), Assertions (4), (5) and (2)). In \cite{13} 
$A$ was an algebra defined over some field $k$ and the
multiplicatively closed set $U\subseteq A[X]$ was the set of $f(X)$ in
$k[X]$ such that $f(0) \neq 0$. Having established Theorem (4.2) the other results of
Section (4) follows, {\it mutatis mutandis}, from \cite{13}. We have
included proofs  in order to make the paper self contained.

\definition{4.1. Notation} We fix  a ring $A$ and a multiplicatively closed subset 
$U\subseteq A[X]$ of the polynomial ring $A[X]$. We write the fraction
ring of $A[X]$ with respect to $U$ as $A[X]_U$.  

Let $ \varphi : A \to
K$ be an $A$-algebra homomorphism. We denote by $U^{\varphi}\subseteq
K[X]$ the image of $U\subseteq A[X]$ under the induced map $A[X] \to
K[X]$. If $f(X)$ is an element in $A[X]$ we denote with
$f^{\varphi}(X)$ the polynomial in $K[X]$ obtained by applying
$\varphi$ to the coefficients of $f(X)$. We have that $U^{\varphi}\subseteq K[X]$ is the set of
polynomials $f^{\varphi}(X)$ where $f(X)$ is in $U\subseteq A[X]$
\enddefinition

\proclaim{Theorem 4.2} Given a monic polynomial $F(X)$ in $A[X]$ which has
positive degree $n$. Let $U\subseteq A[X]$ be a multiplicatively closed
subset. The following three assertions are equivalent.
\roster
\item The canonical map $A[X]/(F(X)) \to A[X]_U/(F(X))$ is an
  isomorphism.
\item The residue classes of $1, X, \ldots ,X^{n-1} $ modulo $F(X)A[X]_U$ form a basis for the
  $A$-module $A[X]_U/(F(X))$.
\item The norm $N_F(f(X))$ with respect to $F(X)$, is a unit in $A$ for
  all $f(X)$ in $U\subseteq A[X]$.
\endroster
\endproclaim

\demo{Proof} It is clear that the two first assertions are
equivalent. We will show that Assertion (1) is equivalent with
Assertion (3). The fraction map $A[X]/(F(X)) \to A[X]_U/(F(X))$ is an
isomorphism if and only if the class of $f(X)$ in $A[X]/(F(X))$ is
invertible for all $f(X)$ in the multiplicatively closed set
$U\subseteq A[X]$. The residue class of $f(X)$ modulo $(F(X))\subseteq A[X]$ is invertible
if and only if the endomorphism $\mu_F (f)$ on $A[X]/(F(X))$ given as
multiplication by the residue class of $f(X)$ is invertible. By Theorem (2.4) we have that
the determinant of the endomorphism $\mu_F(f)$ is the norm $N_F(f(X))$,
and our claim follows.
\enddemo

\proclaim {Corollary 4.3} Let $A$ be  a local ring. Denote the residue field
of $A$ as $K$. Assume that the $K$-vector space
$A[X]_U/(F(X))\otimes_{A}K$ has a basis given by the residue classes of $1, X,
\ldots ,X^{n-1}$. Then we have that the $A$-module $A[X]_U/(F(X))$ has
a basis given by the residue classes of $1, X, \ldots ,X^{n-1}$. In particular
we have that $A[X]_U/(F(X))$ is a finitely generated $A$-module.
\endproclaim

\demo{Proof} Let $\varphi : A \to K$ be the residue class map. We have a canonical isomorphism 
$$ A[X]_U/(F(X))\otimes_A K\cong
K[X]_{U^{\varphi}}/(F^{\varphi}(X)) . \tag{4.3.1}$$
By assumption we have that the residue classes of $1, X, \ldots , X^{n-1}$ is
a $K$-basis for (4.3.1). It then  follows from Assertion (3) of Theorem (4.2) that the norm
$N_{F^{\varphi}}(f^{\varphi}(X)) $ with respect to $F^{\varphi}(X)$,
is a unit in $K$ for all $f^{\varphi}(X)$ in $U^{\varphi}\subseteq
K[X]$. Furthermore we have by Lemma (3.2) that 
$N_{F^{\varphi}}(f^{\varphi}(X))=\varphi (N_F(f(X)))$.
Hence $N_F(f(X))$ is a unit in $A$ for all $f(X)$ in $U\subseteq
A[X]$. Now the claim follows by Assertion (3) of Theorem (4.2).
\enddemo
 
\definition {Remark} It is not true that $A$-modules $A[X]_U/(F(X))$, with
monic polynomials $F(X) \in A[X]$, in
general are finitely generated. For instance, let $A=k[T]$ be the
polynomial ring in the variable $T$ over a field $k$. Let $U\subseteq
A[X]$ be the set of non-zero elements of $A$, thus  $U$ is the set of
non-zero polynomials $f(T)$ in $k[T]$. Let $F(X)=X-T$, which is a monic
polynomial in $A[X]$ of degree 1. We have that $A[X]/(F(X))=A$, from
which it follows that $A[X]_U/(F(X))=k(T)$, the function field of
$A=k[T]$, clearly not finitely generated over $A=k[T]$.
\enddefinition

\proclaim{Lemma 4.4} Let $A=K$ be a field, and let $U\subseteq K[X]$ be
a multiplicatively closed subset. Let $F(X)$ be a polynomial in
$K[X]$. Then we have that the $K$-vector space $K[X]_U/(F(X))$ is
generated by the residue classes of $1, X, \ldots ,X^{n-1}$ modulo
$F(X)K[X]_U$, and  where $n$ is the
degree of $F(X)$.
\endproclaim

\demo{Proof} Denote by $x^i$ the residue class of $X^i$ modulo the
ideal $F(X)K[X]_U$. Taking fraction commutes with taking quotients. Hence
given $f(X)$ in $K[X]$ it is clear that the class of $f(X)$ in
$K[X]_U/(F(X))$ is in the span of $1, x, \ldots, x^{n-1}$. We
must show that given $f(X)$ in $U\subseteq K[X]$ then $f(x)^{-1}$,
the residue class of
$f(X)^{-1}$ modulo $F(X)K[X]_U$,  can be written as a $K$-linear
combination of $1, x, \ldots , x^{n-1}$.

Let $f(X)$ be an element of $U\subseteq K[X]$, which we may assume to
be irreducible. Indeed, if $f(X)$ is
invertible in $K[X]_U$, then it follows that the irreducible factors of $f(X)$ are
invertible in $K[X]_U$. Hence it suffices to  show that the irreducible
factors of $f(X)$ are in the $K$-linear span of $1,x, \ldots
,x^{n-1}$. Consequently we assume that the ideal $(f(X))$ in
$K[X]$ is maximal. 

If $F(X)$ is in the ideal $(f(X))$ then
$F(X)=f(X)G(X)$, where $\text{deg}(G(X))<\text{deg}(F(X))$. We have
that $G(X)$ and $F(X)$ generate the same ideal in $K[X]_U$. Thus by 
replacing $F(X)$ by $G(X)$, we may assume that $F(X)$ is not in the ideal $(f(X))$. That is the
ideals $(F(X))$ and $(f(X))$ are coprime. Hence there exist
polynomials $h(X)$ and $H(X)$ such that
$f(X)h(X)+F(X)H(X)=1$ in $K[X]$. It then follows that 
$$ h(X)+f(X)^{-1}F(X)H(X)=f(X)^{-1}, $$
in $K[X]_U$. From which we get that $h(x)=f(x)^{-1}$ in
$K[X]_U/(F(X))$. The element $h(X)$ is in $K[X]$, hence its residue class
$f(x)^{-1}$ in $K[X]_U/(F(X))$, is in the span of $1, x, \ldots
,x^{n-1}$. We have proven our claim.
\enddemo

\proclaim {Theorem 4.5} Given an ideal $I\subseteq A[X]_U$ such that
the residue ring $A[X]_U/I$ is a free $A$-module of rank $n$. Then
there exists a unique monic polynomial $F(X)$ in $A[X]$ of degree $n$, whose
image in $A[X]_U$ generates the ideal $I$. In particular  we have that the
canonical map
$A[X]/(F(X)) \to A[X]_U/I$
is an isomorphism
\endproclaim

\demo{Proof} Let $V=\oplus_{i=1}^nA$. By assumption
we have that $V=A[X]_U/I$ as $A$-modules. The $A$-module $V$ is also
an $A[X]$-module by the natural
map $ A[X] \to A[X]_U/I$. Thus the variable $X$ is mapped to an
$A$-linear endomorphism $\theta $ on $V$. Let $F(X)=F_{\theta}(X)$ be
the characteristic polynomial of the endomorphism $\theta$.

By the Cayley-Hamilton Theorem (see e.g \cite{11}, XIV \S 3, Theorem
(3.1), p. 561),   any endomorphism satisfies its
characteristic polynomial. Hence we have that $F(X)$ is in the kernel
of the natural map $A[X] \to A[X]_U/I$. We have thus shown that
$$A[X]_U/(F(X)) \to A[X]_U/I \tag {4.5.1}$$
is surjective. We claim that the map (4.5.1) is an isomorphism of $A$-modules.

To show that (4.5.1) is an isomorphism of $A$-modules it is sufficient
to show (4.5.1) when $A$ is a local ring. Let $A$ be a
local ring. Let  $\varphi : A \to K$ be the residue class map. Since the map
in (4.5.1) is surjective it follows that the $K$-vector space
$A[X]_U/(F(X))\otimes_A K $ is of dimension greater or equal to $n$,
the rank of $A[X]_U/I$. We have that $A[X]_U/(F(X))\otimes_A K$ is
isomorphic to $K[X]_{U^{\varphi}}/(F^{\varphi}(X))$. By Lemma (4.4) 
we have that the classes of $1, X, \ldots , X^{n-1}$ generate the $K$-vector space
$A[X]_U/(F(X))\otimes_AK$. Since the dimension of
$A[X]_U/(F(X))\otimes_A K$ is at least $n$ it follows that the
dimension equals $n$. It then
follows from Corollary (4.3) that $A[X]_U/(F(X))$ is a free $A$-module where
the classes of $1, X, \ldots ,X^{n-1}$ form a basis. Thus (4.5.1) is
an isomorphism. 

We have that the classes of $1, X, \ldots ,X^{n-1}$ form a basis for
the $A$-module $A[X]_U/(F(X))$. It follows from Theorem (4.2) that
$A[X]/(F(X))$ is canonically isomorphic to $A[X]_U/(F(X))$. We have
proven the Theorem.
\enddemo

\subhead \S 5. - Parameterizing ideals of fraction rings\endsubhead

In Section (5) we will extract the results from the preceding section
and show that there exist a ring parameterizing the set of ideals in $A[X]_U$ such
that the residue ring is a free, rank $n$ module over $A$.

\definition{5.1. The functor ${\Cal F}^n_U$} We fix a multiplicatively closed subset
$U\subseteq A[X]$. For any $A$-algebra $K$ we let ${\Cal F}^n_U(K)$
denote the set of residue rings $K\otimes_AA[X]_U/I$, which are free
and of rank $n$ as $K$-modules. The map sending an $A$-algebra $K$ to
the set ${\Cal
  F}_U^n(K)$ becomes in a natural way a covariant functor from the category of $A$-algebras to sets. 

The $A$-valued points of ${\Cal F}_U^n$, that is the set
${\Cal F}_U^n(A)$, correspond to  residues of
$A[X]_U$ which are free and of rank $n$ as $A$-modules. We call ${\Cal
  F}_U^n$ the functor parameterizing the free, rank $n$ residues of
$A[X]_U$. The goal of Section (5) is to show that the functor ${\Cal
  F}_U^n$ is representable.
\enddefinition

\definition{5.2. The construction of the universal objects}
To a multiplicatively closed subset $U\subseteq
A[X]$  and a positive integer $n$, we define the
multiplicatively closed subset $U(n)\subseteq \otimes_A^{(n)}A[X]$ as
$$ U(n)=\{ f(X) \otimes \dots \otimes f(X) \mid f(X) \in U \subseteq
A[X]\}.\tag{5.2.1}$$
Let $\xi _U : \otimes_A^{(n)}A[X] \to (\otimes_A^{(n)}A[X])_{U(n)}$
denote the fraction map. Recall (2.3.1) that we have defined
$\Delta_{n,X}(X)=X^n-s_{1,n}(X)X^{n-1}+\dots +(-1)^ns_{n,n}(X)$
in the polynomial ring in the variable $X$ over the ring of symmetric
functions $\otimes _A^{(n)}A[X]$. We denote by
$\Delta_{n,X}^{\xi}(X)$ the element in
$(\otimes_A^{(n)}A[X])_{U(n)}\otimes_AA[X]$ obtained by applying
$\xi_U$ to the coefficients of $\Delta_{n,X}(X)$. Furthermore we define
$$ \aligned V_{n,X}^U &=
V_{n,X}\otimes_{\otimes^{(n)}_AA[X]}(\otimes_A^{(n)}A[X])_{U(n)} \\
&= (\otimes_A^{(n)}A[X])_{U(n)}\otimes_A A[X]/(\Delta_{n,X}^{\xi}(X)).
\endaligned \tag{5.2.2}$$
We have that $V_{n,X}^U$ is a free $(\otimes_A^{(n)}A[X])_{U(n)}$-module
of rank $n$. Note that  $V_{n,X}^{U}$ is not an
$(\otimes_A^{(n)}A[X])_{U(n)}$-valued point of the functor ${\Cal F}^n_U$ since
$V_{n,X}^U$ is not a residue of $(\otimes_A^{(n)}A[X])_{U(n)}\otimes_A
A[X]_U$, but only a residue of $(\otimes_A^{(n)}A[X])_{U(n)}\otimes_A A[X]$.
We need to consider what happens with
$V_{n,X}^U$ when localized, as an $A[X]$-algebra, in the
multiplicatively closed set $U\subseteq A[X]$.
\enddefinition

\proclaim{Lemma 5.3} Let $U\subseteq A[X]$ be a multiplicatively
closed subset. We have that the fraction map of $A[X]$-modules
$$ V_{n,X}^U \to \big(\otimes_A^{(n)}A[X]\big)_{U(n)} \otimes_A
A[X]_U/(\Delta^{\xi_U}_{n,X}(X)),$$
obtained by localization with respect to $U\subseteq A[X]$, is an
isomorphism.
\endproclaim
\demo{Proof} Let $f(X)$  be an element of $A[X]$. Let $\mu (f)$ be the
$\otimes_A^{(n)}A[X]$-linear endomorphism on $V_{n,X}$, given as
  multiplication by the class of $f(X)$ in $V_{n,X}$. It follows by
  Theorem (2.4) with $u_F=\id $ the identity morphism, that the
  determinant of $\mu (f)$ is $f(X)\otimes \dots \otimes f(X)$.  We
  then have  that the determinant of the induced
  endomorphism $\mu (f) \otimes \id$ on $V_{n,X}^U$ is $\xi_U(f(X)
  \otimes \dots \otimes f(X))$. Thus for all $f(X)$ in $U\subseteq
  A[X]$, we have that the determinant of $\mu (f)\otimes \id$ is a
  unit in $(\otimes_A^{(n)}A[X])_{U(n)}$. That is, the class of $f(X)$
  in $V_{n,X}^U$ is invertible, for all $f(X)$ in $U\subseteq
  A[X]$. It follows that the fraction map of the Lemma is an
  isomorphism.
\enddemo

\proclaim{Theorem 5.4} Let $U\subseteq A[X]$ be a multiplicatively
closed subset. For every positive integer $n$ we have that the
functor ${\Cal F}^n_{U}$ parameterizing free, rank $n$ residues of
$A[X]_U$, is represented by the $A$-algebra $\big(
\otimes_A^{(n)}A[X]\big)_{U(n)}$. The universal element is given by
the residue ring  $V_{n,X}^U$.
\endproclaim

\demo{Proof}  We have that
$V_{n,X}^U$ is a free $F=(\otimes_A^{(n)}A[X])_{U(n)}$-module of rank
$n$. It follows  by Lemma (5.3) that $V_{n,X}^U$ is an element of the
set ${\Cal F}^n_{U}(F)$. Consequently we
have a natural transformation $\Phi : \text{Hom}_{A\text{-alg}}(F,-)
\to {\Cal F}^n_{U}$, which, for any $A$-algebra $K$,  sends an $A$-algebra homomorphism
$F \to K$ to the free $K$-module $V_{n,X}^U\otimes_F K$. We
need to show that $\Phi$ is an isomorphism.

Let $ \varphi : A \to K$ be an $A$-algebra homomorphism. Let $I\subseteq K\otimes_AA[X]_U$ be an ideal such such that the residue
ring $K\otimes_AA[X]_U/I$ is a free $K$-module of rank $n$. We have
that $K\otimes_AA[X]_U\cong K[X]_{U^{\varphi}}$, where
$U^{\varphi}\subseteq K[X]$ is the image of $U\subseteq A[X]$ by the
homomorphism $A[X] \to K[X]$. Hence it follows
from Theorem (4.5) that there exist a unique monic polynomial $F(X)$ in
$K[X]$ of degree $n$, such that the image of $F(X)$ in the fraction
ring $K[X]_{U^{\varphi}}$ generates $I$. Moreover, we have that the canonical
fraction map
$$ K[X]/(F(X)) \to K[X]_{U^{\varphi}}/(F(X)) = K\otimes_AA[X]_U/I \tag{5.4.1}$$
is an isomorphism. Let $F(X)=X^n-u_1X^{n-1}+\dots (-1)^nu^n$ in $K[X]$.

We get by  base change an $A$-algebra homomorphism
$\hat \varphi : \otimes_A^{(n)}A[X] \to \otimes_K^{(n)}K[X]$, where
$\varphi : A \to K$ is the structure map. The coefficients of $F(X)$
determine a $K$-algebra homomorphism $u_F : \otimes_K^{(n)}K[X] \to
K$ by $u_F(s_{i,n}(X))=u_i$ for each $i=1, \ldots ,n$. It is clear that the composite morphism $\hat \varphi \circ u_F$
is an $A$-algebra homomorphism such that
$$ V_{n,X}\otimes _{\otimes_A^{(n)}A[X]} K=K[X]/(F(X)). $$
We need to show that the composite map $\hat \varphi \circ u_F$
factors through the fraction ring $F$ of
$\otimes_{A}^{(n)}A[X]$. Since the fraction map of (5.4.1) is an
isomorphism it follows by Assertion (3) of Theorem (4.3), that $N_F(f^{\varphi}(X))$ is
a unit in $K$ for all $f^{\varphi}(X) \in U^{\varphi} \subseteq K[X]$.  It follows that the
symmetric functions $f(X)
\otimes \dots \otimes f(X)$ in $\otimes_A^{(n)}A[X]$ for all $f(X) \in
U\subseteq A[X]$, are mapped to
units in $K$ by the homomorphism $\hat \varphi \circ u_F$. By the
universal property of the fraction ring $F$, the $A$-algebra homomorphism
$\hat \varphi \circ u_F$ factors through $F$.

Thus we have shown that for any $A$-algebra $K$, and  to any element
$V$ in ${\Cal F}^n_U(K)$ there exist an $A$-algebra homomorphism $F
\to K$ such that $V=V_{n,X}^U\otimes_F K$. 

To complete the proof of the Theorem we need
to show that two different elements in $\text{Hom}_{A\text{-alg}}(F,K)$
correspond to two different elements in ${\Cal F}_U^n(K)$. Let
$\rho : F \to K$ be an $A$-algebra homomorphism. Since $F$ is a
fraction ring of $\otimes_A^{(n)}A[X]$, it follows that $\rho$ is
determined by its action on the elementary symmetric functions
$s_{1,n}(X), \ldots , s_{n,n}(X)$. Let $G(X)=X^n-v_1X^{n-1}+\dots
+(-1)^nv_n$, where $v_i=\rho(s_{i,n}(X))$, for each $i=1, \ldots
,n$. We have that that $V^{U}_{n,X}\otimes_F K$ is isomorphic to
$K[X]_{U^{\varphi}}/(G(X))$, which by Theorem (4.5) is isomorphic to
$K[X]/(G(X))$. It is clear that two polynomials $F(X)$ and $G(X)$
in $K[X]$ of degree $n$ and with leading coefficient 1, generate the same ideal if and only if $F(X)=G(X)$. Hence it
follows that two different $A$-algebra homomorphisms $F \to K$ give
two different ideals in $K \otimes_A A[X]_U$. We have proven the Theorem.
\enddemo

\subhead \S 6. - Symmetric products of fraction rings\endsubhead

In the next two sections we study more in detail the  ring
$(\otimes_A^{(n)}A[X])_{U(n)}$ and the ring $V^U_{n,X}$ which together
form the universal pair  of Theorem  (5.4).

In Section (6) we will show that the fraction ring
$(\otimes_A^{(n)}A[X])_{U(n)}$ is canonically isomorphic to the
ring symmetric tensors of $\otimes_A^nA[X]_U$. In Section (7) we show
that $V^U_{n+1,X}$ is isomorphic to
$(\otimes_A^{(n)}A[X])_{U(n)}\otimes_AA[X]_U$. 

\proclaim{Proposition 6.1}Let $K$ and $L$ be two $A$-algebras. Let
$U\subseteq K$ and $V\subseteq L$ be  multiplicatively closed
subsets. Then we have for any $K$-module $M$ and any $L$-module $N$ a
canonical isomorphism of $A$-modules $M_U\otimes_A N_V \cong
(M\otimes_AN)_{U\cdot V}$, where $U\cdot V\subseteq K\otimes_A L$ is the multiplicatively closed
subset given as
$$ U\cdot V =\{f\otimes g \mid f\in U\subseteq K, g \in V \subseteq
L\}.$$
In particular we have that $K_U\otimes_AL_V\cong
(K\otimes_AL)_{U\cdot V}$.
\endproclaim
\demo{Proof} We have by definition that $M_U=M\otimes_K K_U$. Hence
our claim follows if we show that $K_U\otimes_A L_V\cong
(K\otimes_AL)_{U\cdot V}$. 

We have a natural,
well defined map $K\otimes_A L \to K_U\otimes_A L_V$,
sending $\sum_{i=1}^m(f_i\otimes g_i)$ to $\sum_{i=1}^m((f_i,1)\otimes
(g_i,1))$. The map $K\otimes_A  L \to K_U \otimes_A L_V$ maps an element
$(u\otimes v)$ in $U\cdot V$ to a unit in $K_U \otimes_A L_V$. Hence
by the universal property of fraction rings we get a map 
$$ (K\otimes_A L)_{U\cdot V} \to K_U \otimes_A L_V .\tag{6.1.1}$$
The composite map $K \to K\otimes_AL \to (K\otimes_AL)_{U\cdot V}$
sends any element $u$ in $U\subseteq K$ to the unit $u\otimes 1$ in
$(K\otimes_AL )_{U\cdot V}$. Hence we get a natural $A$-algebra homomorphism  $K_U \to
(K\otimes_AL )_{U\cdot V}$. Similarily we get a map $L_V \to
(K\otimes_A L)_{U\cdot V}$. It follows
from the universal property of the tensor product that we have get a
map  $K_U\otimes_A L_V \to (K\otimes_A L)_{U\cdot V}$, easily seen to
be the inverse of (6.1.1).
\enddemo

\proclaim{Lemma 6.2} Let $U$ be a  multiplicatively closed
subset of a ring $A$, containing the identity element. Assume
that $U$ is a subset of a multiplicatively closed set $V$, such that for any element $v$ in $V$ there
exist an element $w$ in $V$ such that $vw$ is in $U$. Then we have that
$A_U$ is canonical isomorphic to $A_V$.
\endproclaim

\demo{Proof} We have by assumption that $U\subseteq V$. Hence there is  
a canonical map $A_U \to A_V$, which is easily seen to be both
injective and surjective.
\enddemo

\definition{6.3. An action of the symmetric group} Let $K$ be an $A$-algebra. The symmetric group of $n$-letters
$\frak S_n$ acts on the tensor algebra $\otimes_A^nK$ in a natural way. If $\sum_{j=1}^m(\otimes_{i=1}^nf_{i,j})$ is
an element of $\otimes_A^nK$ then an element $\sigma $ in the group
$\frak S_n$ acts by 
$\sigma( \sum_{j=1}^m
\otimes_{i=1}^nf_{i,j})=\sum_{j=1}^m(\otimes_{i=1}^nf_{\sigma(i),j})$. The subring of symmetric tensors is written as
$\otimes_A^{(n)}K$.
\enddefinition

\proclaim{Proposition 6.4} Let $K$ be an $A$-algebra, and let $U\subseteq K$ be a multiplicatively
closed set. For every positive integer $n$ we define
$U(n)\subseteq \otimes_A^{(n)}K$ as the multiplicatively closed
set $U(n) =\{f\otimes \cdots \otimes f \mid f \in U \subseteq K\}$.
Then we have a
canonical isomorphism 
$ \otimes _A^{(n)}K_U \cong \big(
\otimes_A^{(n)}K\big)_{U(n)}$. In particular we have that
$$ (\otimes_A^{(n)}A[X])_{U(n)} \cong \otimes_A^{(n)}A[X]_U.$$ 
\endproclaim

\demo{Proof} We define the multiplicatively closed set $U^n\subseteq
\otimes_A^nK$ by 
$$U^n=\{f_1\otimes\cdots \otimes  f_n \mid f_i \in U\subseteq K, \text{ for } i=1, \ldots, n.\}$$ 
It follows by repeated use of  Proposition (6.1) that we have a
canonical isomorphism 
$\otimes_A^nK_U\cong (\otimes_A^nK)_{U^n}$. We have that
$\otimes_A^{(n)}K_U$ is the ring of invariants of $\otimes_A^nK_U$. We will
 show that
the ring of invariants of $(\otimes_A^nK)_{U^n}$ is
$(\otimes_A^{(n)}K)_{U(n)}$. Let $U^{\frak S}=U^n\cap
(\otimes_A^{(n)}K)$. 

Recall (\cite{1}, Exercise 12, p. 68) that if $G$ is a finite group
acting on a ring $B$, such that $G(U)\subseteq U$ for a
multiplicatively closed set $U\subseteq B$. Then we have that
$B^G_{U^G}\cong  (B_U)^G$, where $U^G=U\cap B^G$

Hence we have that the ring of
invariants of $(\otimes_A^{n}K)_{U^n}$
is $(\otimes^{(n)}_AK)_{U^{\frak S}}$.   Clearly we have
$U(n)\subseteq U^{\frak S}$. Our Proposition is proven if we show that
the ring of invariants of $(\otimes_A^nK)_{U^n}$
is isomorphic to $(\otimes_A^{(n)}K)_{U^{\frak S}}$. By Lemma (6.2) it suffices to show that
for any element $f \in $
$U^{\frak S}$, there is an element $H$ in $U^{\frak S}$ such that
the product is in $U(n)$.  Let
$f_{1}\otimes \cdots  \otimes f_{n}$ be  an element of $U^{\frak
  S}\subseteq U^n$. We have that
$$ \prod_{i=1}^n f_i\otimes \cdots \otimes
f_i=\big( f_1\otimes \cdots \otimes f_n\big)H , \tag{6.4.1}$$
for some $H$ in $\otimes_A^nK$. We have that $H$ is in $U^n$, and we
have that
the product (6.4.1) is in $U(n)$. Since the product (6.4.1) is
symmetric and the element $f_1\otimes \dots \otimes f_n$ is
symmetric by assumption,  it follows that $H$ is in $U^{\frak S}$.  We have proven the Proposition. 
\enddemo

\subhead \S 7. The addition map\endsubhead

In Section (7) we show first
that  $\otimes_A^{(n)}A[X]\otimes_A A[X]$ is isomorphic to
$V_{n+1,X}$. Thereafter we show that the isomorphism of $V_{n+1,X}^U$
and $(\otimes_A^{(n)}A[X])_{U(n)}\otimes_AA[X]_U$ follows by
localization.

\definition{7.1. Definition} We have that $\otimes_A^{(n)}A[X]$ is the
polynomial ring over $A$ in the elementary symmetric functions
$s_{1,n}(X), \ldots ,s_{n,n}(X)$ in the variables $X_1, \ldots ,X_n$.  An $A$-algebra
homomorphism from $\otimes_{A}^{(n)}A[X]$  to an $A$-algebra $K$ is
determined by its action on $s_{1,n}(X), \ldots ,s_{n,n}(X)$. For
every positive integer $n$ we 
define the {\it addition} map
$$ a_n :\otimes^{(n)}_AA[X] \to \otimes_A^{(n-1)}A[X] \otimes_A
A[X], \tag{7.1.1}$$
by sending $s_{i,n}(X)$ to $s_{i,n-1}(X)+s_{i-1,n-1}(X)X$ for every
$i=1, \ldots , n$. As a convention we let $s_{0,n}(X)=1$ and
$s_{n,n-1}(X)=0$ for all $n$,
and we set $\otimes_A^{(0)}A[X]=A$. We denote with  
$$\hat a_n : \otimes_A
^{(n)}A[X]\otimes _AA[X] \to  \otimes_A^{(n-1)}A[X]\otimes _AA[X]$$
 the $A[X]$-algebra  induced by the addition map.
\enddefinition

\proclaim{Lemma 7.2} For all positive integers $n$ we define the 
$A[X]$-algebra homomorphism
$$ p_n :\otimes_A^{(n-1)}A[X] \otimes_AA[X] \to 
\otimes_A^{(n)}A[X] \otimes_AA[X]$$
recursively by $p_n(s_{i,n-1}(X))=s_{i,n}(X)-p_n(s_{i-1,n-1}(X))X$ for
$i=1, \ldots ,n$ and
where $p_n(s_{0,n-1}(X))=p_n(1)=1$. Then the following three assertions hold.
\roster
\item We have that the composite map $p_n\circ \hat a_n$ is the
  identity map.
\item We have that $\hat a_n \circ p_n (s_{i,n-1}(X))=s_{i,n-1}(X)$
  for all $i=1, \ldots ,n-1$.
\item The kernel of $\hat a_n$ is generated by $\Delta_{n,X}(X)=\prod_{i=1}^n(X-X_i)$.
\endroster
\endproclaim

\demo{Proof} We first prove Assertion (1). It is enough to show that
$p_n \circ \hat a_n $ is the identity on the elementary symmetric
functions $s_{1,n}(X), \ldots ,s_{n,n}(X)$. For each $i=1, \ldots, n$ we have that 
$$\aligned 
p_n\circ \hat a_n(s_{i,n}(X)) &=p_n
\big(s_{i,n-1}(X)\big)+p_n\big(s_{i-1,n-1}(X) X\big)\\
&=s_{i,n}(X)-p_n(s_{i-1,n-1}(X)X)+p_n(s_{i-1,n-1}(X) X) \\
 &=s_{i,n}(X). 
\endaligned $$
We have proved the first Assertion. The second Assertion is proven by
induction on $i$. For $i=1$ we get by definition that $\hat a_{n}\circ
p_n(s_{1,n-1}(X))=\hat a_{n} (s_{1,n}(X)-X)=s_{1,n-1}(X)$. Assume as
the induction hypothesis that $\hat a_n \circ p_n
(s_{i,n-1}(X))=s_{i,n-1}(X)$ for $i\geq 1$. We then get that
$$\aligned 
\hat a_n \circ p_n(s_{i+1,n-1}(X)) &= \hat a_n \big (s_{i+1,n}(X) -p_n(s_{i,n-1}(X) X) \big) \\
&= s_{i+1,n-1}(X)+s_{i,n-1}(X) X -\hat a_{n}\circ p_n (s_{i,n-1}(X) X)
\\
&=s_{i+1,n-1}(X).
\endaligned $$
Thus we have proven Assertion (2). To prove the last Assertion we
first show that $\Delta_{n,X}(X)$ is in the kernel of $\hat a_n$. We have that
$$ \aligned
\hat a_n(\Delta_{n,X}(X)) &= X^n+\sum_{i=1}^n(-1)^i\hat a_n
(s_{i,n}(X)X^{n-i}) \\ 
&=X^n+\sum_{i=1}^n(-1)^i \big( s_{i,n-1}(X)+s_{i-1,n-1}(X) X \big) X^{n-i} \\
\endaligned  \tag{7.2.1}$$
By definition we have that $s_{0,n-1}(X)=1$ and that
$s_{n,n-1}(X)=0$. Thus it follows from (7.2.1) that $\Delta_{n,X}(X)$
is in the kernel of $\hat a_n$.  By Assertion (1) we have that $\hat
a_n \circ p_n$ is the identity. Since $\Delta_{n,X}(X)$ is in the
kernel of $\hat a_n$ we get an  induced homomorphism
$$ \hat p_n : \big( \otimes_{A}^{(n-1)}A[X]\big ) \otimes_A A[X] \to
V_{n,X}, \tag{7.2.2}$$
where $V_{n,X}=\otimes^{(n)}_AA[X]\otimes_A
A[X]/(\Delta_{n,X}(X))$.  We have that $\hat a_n \circ \hat p_n $ is the
identity map. In particular $\hat a_n$ is surjective. Our claim follows if we prove that the map
$\hat p_n$ is surjective. The $A[X]$-algebra $V_{n,X}$ is
generated by the classes of $s_{1,n}(X), \ldots ,s_{n,n}(X)$. It follows from 
Assertion (2) that we only need to show that the class of  $s_{n,n}(X)$ in
$V_{n,X}$ is in the image of $\hat p_n$.  However we have that 
$$ X^n+s_{1,n}(X)X^{n-1}+\dots +
(-1)^{n-1}s_{n-1,n}(X)X=(-1)^{n-1}s_{n,n}(X)  $$
in $V_{n,X}$. Hence we have that $\hat p_n$ is surjective.  We have
proven the Lemma. 
\enddemo

\proclaim{Lemma 7.3} Let $n$ be a positive integer.  For all integers
$i=1, \ldots, n$ and for all $f(X)$ in $A[X]$ we have that
that the homomorphism $a_n$ sends the element $s_{i,n}(f(X))$ in
$\otimes_A^{(n)}A[X]$ to $s_{i,n-1}(f(X))+s_{i-1,n-1}(f(X))f(X)$ in
$\otimes_A^{(n-1)XS}A[X]\otimes _AA[X]$. In particular we have that $f(X) \otimes
\dots \otimes f(X)$ in $\otimes_A^{(n)}A[X]$ is mapped to  
$$a_n (f(X)\otimes \dots \otimes f(X)) =(f(X)\otimes \dots
\otimes f(X))\otimes f(X)$$.
\endproclaim

\demo{Proof} We have a natural identification $\iota : \otimes_A^nA[X]
\to \otimes_A^{n-1}A[X]\otimes_AA[X]$ which sends $X_i$ to $X_i$ when
$i=1, \ldots , n-1$ and $X_n$ to $X$. For every $i=1, \ldots ,
n$  we have that
$$ \aligned
s_{i,n}(X) & =\sum_{0<k_1< \dots <k_i<n+1}X_{k_1} \cdots
 X_{k_i} \\
&=\sum_{0<k_1<\dots <k_i<n}X_{k_1} \cdots
 X_{k_i}+\sum_{0<k_1<\dots <k_{i-1}<n}X_{k_1} \cdots
X_{k_{i-1}} X_n \\
&=s_{i,n-1}(X)+s_{i-1,n-1}(X) X_n. 
\endaligned \tag{7.3.1}$$
It follows from (7.3.1) that $\iota
(s_{i,n}(f(X)))=s_{i,n-1}(f(X))+s_{i-1,n-1}(f(X))f(X)$. From the
definition of the addition map and (7.3.1) we have that
the restriction of $\iota $ to the subring
$\otimes_A^{(n)}A[X]\subseteq \otimes_A^nA[X]$ coincides with $a_n$,
and our claim follows.
\enddemo

Recall that we defined (5.2.2) the ring $V_{n,X}^U$ as the
localization of the $\otimes_A^{(n)}A[X]$-algebra $V_{n,X}$ with
respect to the multiplicatively closed subset $U(n) \subseteq
\otimes_A^{(n)}A[X]$. We saw in Proposition (6.4) that the fraction
ring $(\otimes_A^{(n)}A[X])_{U(n)}$ was naturally identified with the
ring of symmetric tensors $\otimes_A^{(n)}A[X]_U$. In the next
Proposition we show a similar behavior for $V_{n,X}^U$.

\proclaim{Proposition 7.4} Let $U\subseteq A[X]$ be a multiplicatively
closed subset. We
have that the addition map $a_n : \otimes_A^{(n)}A[X] \to
\otimes_A^{(n-1)}A[X]\otimes_AA[X]$ induces an isomorphism
$$  V_{n,X}^U \cong \big(\otimes_A^{(n-1)}A[X]_U\big)\otimes_AA[X]_U.$$
\endproclaim

\demo{Proof} The homomorphism $\hat a_n :
\otimes_A^{(n)}A[X]\otimes_AA[X] \to \otimes_A^{(n-1)}A[X]\otimes_AA[X]$
is surjective by Assertion (1) of Lemma (7.2). By Assertion (3) of Lemma
(7.2) the kernel of $\hat a_n$ is generated by $\Delta_{n,X}(X)$. Thus
we have that the addition map induces an isomorphism
$$ V_{n,X}\cong (\otimes_A^{(n-1)}A[X])\otimes_AA[X].\tag{7.4.1}$$
When we
localize the $\otimes_A^{(n)}A[X]$-module $V_{n,X}$ with respect to
the multiplicatively closed subset $U(n)\subseteq \otimes_A^{(n)}A[X]$
we get by definition $V_{n,X}^U$. Consequently the proof of
the Proposition will be complete when we show that tensoring the right
term in (7.4.1)
with $\big(\otimes_A^{(n)}A[X]\big)_{U(n)}$ gives 
$(\otimes_A^{(n-1)}A[X]_U\big)\otimes_AA[X]_U$.

Let $W\subseteq B=\otimes_A^{(n-1)}A[X]\otimes_AA[X]$ be the
multiplicatively closed set
$$W=\{(f(X)\otimes \dots \otimes f(X))\otimes f(X) \mid f(X) \in
U\subseteq A[X]\}.$$
It follows from Lemma (7.3) that $W$ is the image of the
multiplicatively closed subset $U(n)\subseteq \otimes_A^{(n)}A[X]$ by
the addition map $a_n : \otimes_A^{(n)}A[X] \to B$. Hence we have that
$B\otimes_{\otimes_A^{(n)}A[X]}\big( \otimes_A^{(n)}A[X]\big)_{U(n)}=B_W$. Furthermore we have by Proposition
(6.4) that $\otimes_A^{(n-1)}A[X]_U$ is the fraction ring of
$\otimes_A^{(n-1)}A[X]$ with respect to the multiplicatively closed
subset $U(n-1)$. It then follows from Proposition (6.1) that we have
$$ (\otimes_A^{(n-1)}A[X]_U)\otimes_A A[X]_U\cong
(\otimes_A^{(n-1)}A[X]\otimes_AA[X])_{U(n-1)\cdot U},$$
where $U(n-1)\cdot U=\{F(X)\otimes g(X) \mid F(X) \in U(n-1), g(X) \in
U\}$. It is clear that we have $W\subseteq U(n-1)\cdot U$. Hence to
prove the Proposition it is
by Lemma (6.2) sufficient to show that for any $f$ in $U(n-1)\cdot U$
there is $g$ in $U(n-1)\cdot U$ such that the product is contained
in $W$.

Given $F(X)\otimes g(X)$ in $U(n-1)\cdot U$, where $F(X)=f(X)\otimes
\dots \otimes f(X)$ is in $U(n-1)$, and where $g(X)$ is in $U$. In the
ring $B$ we have the element $g(X)\otimes \dots \otimes g(X) \otimes
f(X)$, easily seen to be in $U(n-1)\cdot U$. We have that
$$ \aligned \Big (F(X)\otimes g(X)\Big) \cdot \big( g(X)\otimes \dots
\otimes g(X) \otimes f(X) \big)
\\
f(X)g(X) \otimes \dots \otimes f(X)g(X),\endaligned
$$
which is an element of $W$.  We have proven the Proposition. 
\enddemo

\definition{Remark} B. Iversen defines (See \cite{9}, page 3, Section
(1.4)) for any flat $A$-module $C$ a
canonical map 
$$\otimes_A^{(n+m)}C \to \big(\otimes_A^{(n)}C\big)\otimes_A
\big(\otimes_A^{(m)}C\big), $$
which sends $(x_1\otimes \dots \otimes x_{n+m})$ to $(x_1\otimes \dots
\otimes x_m)\otimes (x_{m+1}\otimes \dots \otimes x_{n+m})$. When
applied to our situation with $C=A[X]$, and $m=1$, we have a map
$\otimes_A^{(n+1)}A[X] \to \otimes_A^{(n)}A[X] \otimes _A A[X]$. It follows from Lemma (7.3) that when applied to our situation 
then  our addition map $ a_{n+1}$ coincides with the canonical
map of B. Iversen.
\enddefinition

\subhead \S 8. - Application to Hilbert schemes\endsubhead

\definition{8.1. The Hilbert functor of points} Let $A$ be a
commutative ring. For any $A$-algebra $ A \to K$ we define the
 Hilbert functor ${\Cal Hilb}_{K/A}^n$ of $n$-points on
$\text{Spec}(K)$ as the contravariant functor from the category of
$A$-schemes to sets, mapping an $A$-scheme $T\to \text{Spec}(A)$ to the set

$$ \aligned
{\Cal Hilb}^n_{K/A}(T) \endaligned =  \left\{ \aligned & \text{Closed
  subschemes } Z\subseteq T\times_A \text{Spec}(K) \text{ such that}  \\
&\text{the projection map } p: Z \to T \text{ is flat and where}\\
&\text{the global sections of }p^{-1}(t) \text{ is of dimension }
n\text{ as}\\ 
&\text{a }\kappa (t)\text{-vector space, for all points } t\in T.
\endaligned \right\} $$

We call ${\Cal Hilb}^n_{K/A}$ the
{\it Hilbert functor} of $n$-points on $\text{Spec}(K)$. In those
cases when the functor ${\Cal Hilb}^n_{K/A}$ is representable we call
the representing scheme the {\it Hilbert scheme} of $n$-points on $\text{Spec}(K)$.

If $K$ is an $A$-algebra, we have the $A$-algebra of symmetric tensors
$\otimes_A^{(n)}K$. We define the $A$-scheme 
$$ \text{Symm}_A^n(K)=\text{Spec}(\otimes_A^{(n)}K). $$
\enddefinition

\proclaim{Theorem 8.2} Let 
$U\subseteq A[X]$ be a multiplicatively closed subset of the
polynomial ring in the variable $X$ over a ring $A$. Let $n$ be a fixed positive integer. Then we have
that the $A$-scheme $\text{Symm}^n_{A}(A[X]_U) \to
\text{Spec}(A)$ represents the functor ${\Cal
  Hilb}^n_{A[X]_U/A}$. The universal family is given as
$$\text{Symm}_A^{n-1}(A[X]_U)\times_A \text{Spec}(A[X]_U) \to \text{Symm}^n_A(A[X]_U).$$
\endproclaim

\demo{Proof} First we show  that
$\text{Symm}_A^{n-1}(A[X]_U)\times_A \text{Spec}(A[X]_U)$ is an
element of ${\Cal Hilb}^n_{A[X]_U/A}(\text{Symm}^n_A(A[X]_U))$. 
By Proposition (7.4) we have that the coordinate ring of
$\text{Symm}_A^{n-1}(A[X]_U)\times_A \text{Spec}(A[X]_U)$ is
isomorphic to $V_{n,X}^U$. As a consequence of Theorem (5.3) we have
that $V_{n,X}^U$ is a free, rank $n$ residue of
$(\otimes_A^{(n)}A[X])_{U(n)}\otimes_AA[X]_U$. Finally we have by Proposition
(6.4) that $(\otimes_A^{(n)}A[X])_{U(n)}\cong
\otimes_A^{(n)}A[X]_U$.  Thus we have that
$\text{Symm}_A^{(n-1)}(A[X]_U)\times_A \text{Spec}(A[X]_U)$ is a
$\text{Symm}^n_A(A[X]_U)$-valued point of the Hilbert functor of
$n$-points on $A[X]_U$. 

The next step in the proof is to show that the induced natural
transformation of functors $\text{Hom}(-,\text{Symm}^n_A (A[X]_U)) \to
{\Cal Hilb}^n_{A[X]_U/A}$ is an isomorphism. Given a $A$-scheme $T$
and let $Z$ be a $T$-valued point of ${\Cal Hilb}^n_{A[X]_U/A}$.
Let $\text{Spec}(K)\subseteq T$ be an open affine subscheme. Let
$\varphi : A \to K$ be the $A$-algebra homomorphism corresponding to
the structure map $\text{Spec}(K) \to \text{Spec}(A)$. Let the inverse
image  
$Z\times_{T}\text{Spec}(K)$ be given by the ideal $I\subseteq K\otimes_A
A[X]_U$. 

We have that $K\otimes_A A[X]_U \cong K[X]_{U^{\varphi}}$, where
$U^{\varphi}\subseteq K[X]$ is the image of the multiplicatively
closed set $U\subseteq A[X]$ under the induced map $A[X] \to
K[X]$. Hence the $K$-algebra $K\otimes_A A[X]_U$  is essentially of finite type. By assumption we have that the $K$-module $M=K\otimes_A
A[X]_U/I$ is a flat $K$-module such that for each prime ideal $P$ in
$K$ we have that the $\kappa (P)=K_P/PK_P$-vector space  $M\otimes
\kappa (P)$ is of dimension $n$. It follows (see \cite{12}, Theorem (3.5)) that $M$ is locally a free $K$-module. By
possibly shrinking the open set $\text{Spec}(K)\subseteq T$, we may
assume that $M$ is a free $K$-module of rank $n$.

We have that there exist a covering $\{U_i\}_{i\in {\Cal I}}$ of $T$,
where $U_i\subseteq T$ is open and affine, such that for every $i\in
{\Cal I}$ we have that $p_{i*}({\Cal O}_{Z\times_TU_i})$ is a free
${\Cal O}_{U_i}$-module of rank $n$. Here $p_i: U_i\times_TZ \to
U_i$ is the projection map. It  follows by Theorem (6.3) that there exist a unique
$A$-morphism $f_i: U_i \to \text{Symm}^n_A(A[X]_U)$ such that
$$ U_i\times_TZ \cong
U_i\times_{\text{Symm}^n_A(A[X]_U)}\text{Symm}_A^{n-1}(A[X]_U)\times_A
\text{Spec}(A[X]_U),\tag{8.2.1}$$
for every $i\in {\Cal I}$. It follows from the uniqueness of the
morphisms $f_i$ that they glue together and give a unique morphism $f_Z:
T \to \text{Symm}^n_A(A[X]_U)$ such that the element $Z \in {\Cal
  Hilb}^n_{A[X]_U/A}(T)$ is the pull-back of
$\text{Symm}^{(n-1)}_A(A[X]_U)\times_A \text{Spec}(A[X]_U)$ by the
morphism $f_Z$. We have  proven the Theorem.
\enddemo

\definition{Remark} Given a projective morphism  $C\to S$, where we fix an
embedding of $C$ in some projective $N$-space $\bold P^N_S$ over
$S$. Consider now only $S$-schemes that are locally noetherian. From
A. Grothendiecks general theory of Hilbert functors, we have that
${\Cal Hilb}^n_{C/S}$ is representable (and in fact projective) \cite{5}. In the special case when $C \to S$ is
projective, smooth and of relative dimension 1 over
$S$, it was remarked (see \cite{5}, p. 275)  that the Hilbert
functor ${\Cal Hilb}^n_{C/S}$ is represented by
$\text{Symm}^n_{S}(C)$. 
\enddefinition

\definition{Remark. A comparison with B. Iversens theory of $n$-fold sections} Let $C \to S$ be a flat morphism of schemes. If $T \to S$
is a morphism of schemes then an $n$-fold section of $C$ over $T$ is a
closed subscheme $Z\subseteq T\times_SC$ such that the projection map
$Z \to T$ is finite, flat and of rank $n$. We denote with ${\Cal
  F}^n_{C/S}(T)$ the set of $n$-fold sections of $C$ over $T$. It is
clear that we have a contravariant functor ${\Cal F}^n_{C/S}$ from the
category of $S$-schemes to sets, \cite{9}.

We will compare the functor ${\Cal F}^n_{C/S}$ with the Hilbert
functor ${\Cal Hilb}^n_{K/A}$ in the specific situation when
$C=\text{Spec}(K)$ and $S=\text{Spec}(A)$. It is clear that ${\Cal
F}^n_{\Spec(K)/\Spec(A)}$ is a subfunctor of ${\Cal Hilb}^n_{K/A}$.

I will show, by repeating an argument given in the proof of the
Theorem, that when $K$ is a flat $A$-algebra, essentially of finite
type, then the two functors ${\Cal F}^n_{\Spec(K)/\Spec(A)}$ and
${\Cal Hilb}^n_{K/A}$ are naturally identified. Note that if $Z$ is a
$T$-valued point of ${\Cal Hilb}^n_{K/A}$ then it is not obvious that the
projection map $Z \to T$ is finite.

Let $\text{Spec}(R)$ be an open
subscheme of an $A$-scheme $T$.
Let $Z$ be an $T$-valued point of ${\Cal Hilb}^n_{K/A}$.  We have that $Z$
is a closed subscheme of $T\times_A \text{Spec}(K)$, hence $Z_R =Z\times
_T \text{Spec}(R)$ is a closed subscheme of $\text{Spec}(R)\times_A
\Spec (K)$. Let $Z_R=\text{Spec}(M)$. It follows from the definition of the Hilbert functor
${\Cal Hilb}^n_{K/A}$ that $M$ is a flat $R$-module such that for all
prime ideal $P$ in $R$ we have that the $\kappa (P)=R_P/PR_P$-vector
space $M\otimes_R\kappa (P)$ is of dimension $n$. Since $K$ is
essentially of finite type over $A$, it follows that $M$ is
essentially of finite type over $R$. Hence it follows (\cite{12},
Theorem (3.5)), that $M$ is locally free over $R$. The rank of $M$ as
an $R$-module is clearly $n$, and consequently the $T$-valued point
$Z$ of ${\Cal Hilb}^n_{K/A}$ is a $T$-valued point of ${\Cal F}^n_{\Spec(K)/\Spec(A)}$.

We thus have two functors ${\Cal F}^n_{C/S}$ and ${\Cal Hilb}^n_{K/A}$
that are equal when $K$ is a flat $A$-algebra, essentially of finite
type. In particular the functors are equal when $K=A[X]_U$, a fraction
ring of the polynomial ring in the variable $X$ over $A$. It is
therefore natural to compare the results of \cite{9} with the ones
have in the present paper. In particular there are two results that
are closely related to the Main Theorem of the present paper.

One of the results which B. Iversen shows (\cite{9}, Theorem (3.4),
p. 26) is the following. 

Let $C \to S$ be a flat and finite morphism of schemes such
that for all points $s\in S$, any finite set of points of the fiber over
$s$ is contained in an open affine subset of $C$, whose image by $C\to
S$ is contained in an open affine of $S$ (\cite{9}, page 21, Section
(1.1)). If the canonical map 
$\text{Symm}^{(n-1)}_S(C)\times_S C \to \text{Symm}^{(n)}_S(C)$ is finite,
flat and of rank $n$, then the functor ${\Cal F}^n_{C/S}$ is
represented by $\text{Symm}^{(n)}_S(C)$ and where the universal family is
$\text{Symm}^{(n-1)}_S(C)\times_SC$. The canonical map is the one we discussed in the
Remark of Section (7) in the present paper.

Hence, modulo the result (Theorem (3.4)) of B. Iversen, we could have shortened our
proof of Theorem (8.2). Because to prove the statement of Theorem
(8.2) it would  be sufficient to show that
$\text{Symm}_A^n(A[X]_U)\times_A \text{Spec}(A[X]_U)$ is an
$\text{Symm}^n_A(A[X]_U)$-valued point of
${\Cal Hilb}^n_{A[X]_U/A}$. 

The other result of \cite{9} (Proposition (4.1), p. 29), which I want
to mention is the case when $C\to S$ is a flat family of smooth curves
(and where $S$ is locally noetherian). With these extra assumptions
B. Iversen shows that the canonical map
$\text{Symm}^{(n-1)}_S(C)\times_S C \to \text{Symm}^{(n)}_S(C)$ is
finite, flat and of rank $n$. In other words the 
functor ${\Cal F}^n_{C/S}$ is represented by
$\text{Symm}^{(n)}_S(C)$. 

Comparing with our situation we have that
the homomorphism $\Spec(A[X]_U) \to \Spec(A)$ is smooth,
but the fibers are not necessarily of finite type, hence neither curves. When $S=\text{Spec}(A)$ is noetherian
and $A[X]_U$ is finitely generated as an $A$-algebra, or equivalently
that $\Spec (A[X]_U)$ is an basic open subset of $\Spec(A[X])$, then our Theorem
(8.2) is a consequence of \cite{9} (Proposition (4.1), p. 29). 
\enddefinition

\definition{8.3. Example I: The Hilbert scheme of points on the affine
  line} When $U=\{1\}\subseteq A[X]$ is the trivial subset we have that
  $\text{Spec}(A[X]_U)=\A^1_A$ is the affine line over
  $A$.  The Hilbert scheme of $n$-points on $\A^1_A$ is
  given as $\text{Spec}(\otimes_A^{(n)}A[X])$. The ring of symmetric
  functions $\otimes_A^{(n)}A[X]$ is the polynomial ring in
  $n$-variables over $A$. Thus the parameterizing scheme $\text{Symm}^n_A(A[X])$ is simply the
  affine $n$-space $\A^n_A$ over $A$. Note that the only assumptions
  on the base ring $A$ is that $A$ is commutative and unitary.
\enddefinition

\definition{8.4. Example II: The Hilbert scheme of points on open
  subsets of the line} Let the multiplicatively closed subset $U$ be
  given by multiples of one element $f$ in $A[X]$, that is
  $U=\{f^m\}_{m\geq 0}$. Then $\text{Spec}(A[X]_U)$ is an basic open
  (possibly empty) subscheme of $\A^1_A$, the affine line over $A$. 

The Hilbert scheme of $n$-points on $\text{Spec}(A[X]_U)$ is given as
the spectrum of $(\otimes_A^{(n)}A[X])_{U(n)}$, where $U(n)=\{ (f\otimes
\dots \otimes f)^m\}_{m\geq 0}$. Hence we have that the Hilbert functor of
$n$-points on an basic open subscheme of the line is represented by an
open subscheme of the Hilbert scheme of $n$-points on the line.
\enddefinition

\definition{8.5. Example III: The Hilbert scheme parameterizing finite
  length subschemes of the line with support at the origin} Let the base ring $A=k$ be a field, and
  let $U\subseteq k[X]$ be the set of polynomials $f(X)$ such that
  $f(0)\neq 0$. Thus $k[X]_U=k[X]_{(X)}$ is the local ring of the
  origin on the line, and the functor ${\Cal
  Hilb}^n_{k[X]_{(X)}/k}$ parameterizes the length $n$ subschemes of
  $\text{Spec}(k[X]_{(X)})$. There is only one closed subscheme of $\text{Spec}(k[X]_{(X)})$ of
length $n$, namely the scheme given by the ideal $(X^n)\subseteq
k[X]_{(X)}$. A situation which was studied in detail in \cite{13}.

\enddefinition

\definition{8.6. Example IV: A Hilbert scheme without rational points}
Assume that the base ring $A$ is an integral domain. Let $U\subseteq
A[X]$ be the set of non-zero polynomials. We have that $A[X]_U=A(X)$,
the field of fractions of $A[X]$. We have that $\text{Symm}^n_A(A(X))$
represents the Hilbert functor of $n$-points on
$\text{Spec}(A(X))$. Note that $\text{Spec}(A(X))$ is just a
point and that there exist no non-trivial subschemes of
$\text{Spec}(A(X))$. Consequently the Hilbert scheme of $n$-points on
$\text{Spec}(A(X))$ has no $A$-valued points. In particular if $A=k$
is a field, we have that the Hilbert scheme of $n$-points on
$\text{Spec}(k(X))$ has no $k$-rational points.

We will show that the parameterizing
scheme $\text{Symm}_k^n(k(X))$ is of dimension $(n-1)$ when the base
ring $A=k$ is a field. 

We have that the coordinate ring is, Proposition (6.4),
$(\otimes^{(n)}_kk[X])_{U(n)}$, where $U(n)$ is the set of products of
the form
$f(X)\otimes \dots \otimes f(X)$, with non-zero polynomials $f(X) \in
k[X]$. 

The ring of symmetric functions $\otimes_k^{(n)}k[X]$ is the
polynomial ring in the variables $s_{1,n}(X), \ldots ,s_{n,n}(X)$, and
is consequently of
dimension $n$. It is clear that to any maximal ideal $P$ in
$\otimes_k^{(n)}k[X]$, we can find an element of the form
$f(X)\otimes \dots \otimes f(X)$, with $f(X)\neq 0$, in $P$. Thus the
extension of the maximal ideals in $\otimes_k^{(n)}k[X]$, in
the fraction ring $(\otimes^{(n)}_kk[X])_{U(n)}$ becomes the whole
ring. A phenomena which we expected since the parameterizing scheme
$\text{Symm}^n_k(k(X))$ has no $k$-rational points.

We have that all the maximal ideals of $\otimes_k^{(n)}k[X]$ meets the set
$U(n)$. It follows that we  the dimension of
$(\otimes^{(n)}_kk[X])_{U(n)}$ is at most $n-1$. Next we note that
if $G(X_1, \ldots ,X_n)$ is an irreducible  polynomial in the variables
$X_1, \ldots ,X_n$, which is symmetric, then $G(X_1,\ldots ,X_n)$
correspond to a height 1 prime ideal in $\otimes_k^{(n)}k[X]$ which
does not meet $U(n)$. If $n>1$ clearly such functions exist. The elementary
symmetric functions $s_{1,n}(X), \ldots ,s_{n-1,n}(X)$ are examples of such.
Thus we have that the ideal $P$ generated by $s_{1,n}(X), \ldots ,s_{n-1,n}(X)$
does not meet $U(n)$. Consequently the extension of $P$ in the
fraction ring $(\otimes_k^{(n)}k[X])_{U(n)}$ correspond to a prime
ideal. We have that the localization of
$\otimes_k^{(n)}k[X]$ in $P$ is a local ring of dimension
$(n-1)$. Since $(n-1)$ was an upper bound for the dimension of the
fraction ring $(\otimes_k^{(n)}k[X])_{U(n)}$, it follows that $(n-1)$
is the dimension of $\text{Symm}_k^n(k(X))$. 

It may be surprising that we need a $(n-1)$ dimensional scheme to
parameterize the empty set of closed subschemes of $\text{Spec}(k(X))$
of finite length.
\enddefinition

\definition{8.7. Example V: Hilbert schemes of one point} Let the fixed
integer $n=1$. For any multiplicatively closed subset $U\subseteq
A[X]$ we have that
$\text{Symm}^{1}_A(A[X]_U)=\text{Spec}(A[X]_U)$. Hence the scheme
$\text{Spec}(A[X]_U)$ itself is the Hilbert scheme of one point on
$\text{Spec}(A[X]_U)$. See also  \cite{10} (Corollary (2.3) of
Proposition (2.2), p. 109) where S.L.  Kleiman proves that for any $S$-scheme
$X$ the functor ${\Cal
  Hilb}^1_{X/S}$ is represented by the scheme $X$, and where the
universal family is given by the diagonal in $X\times_S X$. 
\enddefinition

\enddocument
\end